\documentclass{amsart}
\usepackage{amssymb, amsmath, latexsym}

\renewcommand{\baselinestretch}{\baselinestretch}
\renewcommand{\baselinestretch}{1.1}
\author{Constantin-Nicolae Beli}
\title[]{On the kernel of the projection map $T(V)\to S(V)$} 
\date{}

   \def\m{\lim}

 \def\({\overline}

\def\){\underline} \def\<{\cdot} 
\def\>{~~~~~~~} \def\#{{\bf
Definition}} \def\*{\section} \def\be{\begin{equation}}
\def\ee{\end{equation}}

  \def\sbq{\subseteq} \def\spq{\supseteq}

\def\bmat{\left(\begin{array}} \def\emat{\end{array}\right)}

 \def\m2{~(\mo 2)} \def\no{\noindent}
 \def\btm{\begin{thm}}
\def\etm{\end{tm}}
 \def\blem{\begin{lem}}
\def\elem{\end{lem}}

\newtheorem{theorem}{Theorem}[section]
\newtheorem{proposition}[theorem]{Proposition}
\newtheorem{lemma}[theorem]{Lemma}
\newtheorem{definition}{Definition}
\newtheorem{corollary}[theorem]{Corollary}

\newtheorem{bof}[theorem]{}
\newtheorem{teorema}{Theorem}

\def\qed{\mbox{$\Box$}\vspace{\baselineskip}}
\def\pf{$Proof.$ } 
\def\bco{\begin{corollary}} \def\eco{\end{corollary}} 
\def\bdf{\begin{definition}} \def\edf{\end{definition}} 
\def\btm{\begin{theorem}} \def\etm{\end{theorem}} 
\def\bpr{\begin{proposition}} \def\epr{\end{proposition}}  
\def\blm{\begin{lemma}} \def\elm{\end{lemma}} 
\def\bff{\begin{bof}\rm} \def\eff{\end{bof}}
\def\btr{\begin{teorema}} \def\etr{\end{teorema}}

\def\de{\newcommand} \de\tm[1]{{\no\bf Theorem~#1}} 
\def\mb{\mathbb} 
   \def\ZZ{{\mb Z}}
 
\de\lm[1]{{\no\bf Lemma~#1}}
\de\df[1]{{\no\bf Definition~#1}} \de\co[1]{{\no\bf Corollary~#1}}
\de\lr[1]{\longrightarrow^{\!\!\!\!\!\!\!\! #1}}
\de\lf[1]{\longleftarrow^{\!\!\!\!\!\!\!\! #1}}
\de\si[1]{\sim^{\!\!\!\!\! #1}} \de\apr[1]{\approx^{\!\!\!\!\! #1}}
\de\leg[2]{\left(\frac {#1}{#2}\right)}
\DeclareMathOperator\Br{Br}

\de\Brr[1]{{}_{#1}\Br}

\DeclareMathOperator\Ima{Im}

\begin{document}

\begin{abstract}

Let $V$ be a vector space over some field $F$ and let
$\rho_{T,S}:T(V)\to S(V)$ be the projection map, given by
$x_1\otimes\cdots\otimes x_n\mapsto x_1\cdots x_n$.

In this paper we give a descrption of $\ker\rho_{S,T}$ in terms of
generators and relations. Namely, we will define a $\ZZ_{\geq
2}$-graded $T(V)$-bimodule $M(V)$, which is a quotient of the
$T(V)$-bimodule $T(V)\otimes\Lambda^2(V)\otimes T(V)$, and a morphism
of $T(V)$-bimodules $\rho_{M,T}:M(V)\to T(V)$, such that the sequence
$$0\to M(V)\xrightarrow{\rho_{M,T}}T(V)\xrightarrow{\rho_{T,S}}S(V)\to
0$$
is exact.

In a related result, we define the algebra $S'(V)$ as $T(V)$
factorised by the bilateral ideal generated by $x\otimes y\otimes
z-y\otimes z\otimes x$, with $x,y,z\in V$, and we prove that there is a
short exact sequence,
$$0\to\Lambda^{\geq 2}(V)\xrightarrow{\rho_{\Lambda^{\geq
2},S'}}S'(V)\xrightarrow{\rho_{S',S}}S(V)\to 0.$$

When considering the homogeneous components of degree $2$, we have
$M^2(V)=\Lambda^2(V)$ and $S'^2(V)=T^2(V)$ so in both cases we get the
well known exact sequence
$$0\to\Lambda^2(V)\to T^2(V)\to S^2(V)\to 0.$$
\end{abstract}
\maketitle

\section{The bimodule $M(V)$}
Let $V$ a vector space over a field $F$ and let $(v_i)_{i\in I}$ be a
basis, where $(I,\leq )$ is a totally ordered set. 

Let $\rho_{T,S}:T(V)\to S(V)$ be the canonical projection and, for
$n\geq 0$, let $\rho_{T^n,S^n}:T^n(V)\to S^n(V)$ be its homogeneous
component of degree $n$. 

We denote by $[\cdot,\cdot ]:T(V)\times T(V)\to T(V)$ the commutator
map, $[\xi,\eta ]=\xi\otimes\eta -\eta\otimes\xi$. Then
$\ker\rho_{T,S}$ is the bilateral ideal generated by $[x,y]$, with
$x,y\in V$. We want to describe $\ker\rho_{T,S}$ in terms of
generators and relations. On homogeneous components, for $n=0,1$ we
have $T^n(V)=S^n(V)$ so $\ker\rho_{T^n,S^n}=0$. The first interesting
case is $n=2$. The map $V^2\to T^2(V)$, given by $(x,y)\mapsto [x,y]$,
is bilinear and alternating, so it induces a linear map
$\rho_{\Lambda^2,T^2}:\Lambda^2(V)\to T^2(V)$, given by $x\wedge
y\mapsto [x,y]$. Then we have the following well known and elementary
result.

\bpr We have an exact sequence
$$0\to\Lambda^2(V)\xrightarrow{\rho_{\Lambda^2,T^2}}
T^2(V)\xrightarrow{\rho_{T^2,S^2}} S^2(V)\to 0.$$
\epr

We now consider the $T(V)$-bimodule $T(V)\otimes\Lambda^2(V)\otimes
T(V)$, generated by $\Lambda^2(V)$. It is $\ZZ_{\geq 2}$ graded,
where for every $n\geq 2$ homogeneous component of degree $n$ is
$$(T(V)\otimes\Lambda^2(V)\otimes
T(V))^n=\bigoplus_{i+j=n-2}T^i(V)\otimes\Lambda^2(V)\otimes T^j(V).$$

Note that we can define the commutator $[\cdot,\cdot ]$, by the same
formula, $[\xi,\eta ]=\xi\otimes\eta -\eta\otimes\xi$, also as
$[\cdot,\cdot ]:T(V)\times (T(V)\otimes\Lambda^2(V)\otimes T(V))\to
T(V)\otimes\Lambda^2(V)\otimes T(V)$, or as $[\cdot,\cdot
]:(T(V)\otimes\Lambda^2(V)\otimes T(V))\times T(V)\to
T(V)\otimes\Lambda^2(V)\otimes T(V)$.

We consider the map $1\otimes\rho_{\Lambda^2,T^2}\otimes
1:T(V)\otimes\Lambda^2(V)\otimes T(V)\to T(V)$. Note that
$1\otimes\rho_{\Lambda^2,T^2}\otimes 1$ is a morphism of graded
$T(V)$-bimodules.

Now $T(V)\otimes\Lambda^2(V)\otimes T(V)$ is spanned by $\xi\otimes
x\wedge y\otimes\eta$, with $x,y\in V$ and $\xi,\eta\in T(V)$, so
$\Ima (1\otimes\rho_{\Lambda^2,T^2}\otimes 1)$ is spanned by
$(1\otimes\rho_{\Lambda^2,T^2}\otimes 1)(\xi\otimes x\wedge
y\otimes\eta)=\xi\otimes [x,y]\otimes\eta$, i.e. it is the ideal of
$T(V)$ generated by $[x,y]$, with $x,y\in V$. Hence $\Ima
(1\otimes\rho_{\Lambda^2,T^2}\otimes 1)=\ker\rho_{T,S}$ and we have
the exact sequence
$$T(V)\otimes\Lambda^2(V)\otimes
T(V)\xrightarrow{1\otimes\rho_{\Lambda^2,T^2}\otimes
1}T(V)\xrightarrow{\rho_{T,S}}S(V)\to 0.$$

It follows that $\ker\rho_{T,S}\cong\frac{T(V)\otimes\Lambda^2(V)\otimes
T(V)}{\ker (1\otimes\rho_{\Lambda^2,T^2}\otimes 1)}$.

\blm The following elements of $T(V)\otimes\Lambda^2(V)\otimes T(V)$
belong to\\ $\ker (1\otimes\rho_{\Lambda^2,T^2}\otimes 1)$.

(i) $[x,y]\otimes\xi\otimes z\wedge t-x\wedge y\otimes\xi\otimes
[z,t]$, with $x,y,z,t\in V$, $\xi\in T(V)$.

(ii) $[x,y\wedge z]+[y,z\wedge x]+[z,x\wedge y]$, with $x,y,z\in V$.
\elm

\pf (i) We have
\begin{multline*}
(1\otimes\rho_{\Lambda^2,T^2}\otimes 1)([x,y]\otimes\xi\otimes z\wedge
t-x\wedge y\otimes\xi\otimes [z,t])\\
=[x,y]\otimes\xi\otimes [z,t]-[x,y]\otimes\xi\otimes [z,t]=0.
\end{multline*}

(ii) By the Jacobi identity we have
$$(1\otimes\rho_{\Lambda^2,T^2}\otimes 1)([x,y\wedge z]+[y,z\wedge
x]+[z,x\wedge y])=[x,[y,z]]+[y,[z,x]]+[z,[x,y]]=0.$$
\qed

\bdf Let $M(V)=(T(V)\otimes\Lambda^2(V)\otimes T(V))/W_M(V)$, where
$W_M(V)$ is the subbimodule of $T(V)\otimes\Lambda^2(V)\otimes T(V)$
generated by $[x,y]\otimes\xi\otimes z\wedge t-x\wedge
y\otimes\xi\otimes [z,t]$, with $x,y,z,t\in V$ and $\xi\in T(V)$, and
$[x,y\wedge z]+[y,z\wedge x]+[z,x\wedge y]$, with $x,y,z\in V$.

If $\eta\in T(V)\otimes\Lambda^2(V)\otimes T(V)$ then we denote by
$[\eta ]$ its class in $M(V)$.

On $M(V)$ we keep the notation $\otimes$ for the left and right multiplication
from $T(V)\otimes\Lambda^2(V)\otimes T(V)$. That is $\xi\otimes
[\eta ]\otimes\xi':=[\xi\otimes\eta\otimes\xi']$ $\forall\xi,\xi'\in
T(V)$, $[\eta ]\in M(V)$.
\edf

Note that $W_M(V)$ is generated by homogeneous elements so it is
homogeneous. (In the formula $[x,y]\otimes\xi\otimes z\wedge t-x\wedge
y\otimes\xi\otimes [z,t]$ we may restrict ourselves to $\xi\in
T^k(V)$, with $k\geq 0$, which makes it homogeneous of degree $k+4$.)
Therefore $M(V)$ inherits from $T(V)\otimes\Lambda^2(V)\otimes T(V)$
the property of being a $\ZZ_{\geq 2}$-graded $T(V)$-bimodule.

By Lemma 1.2, we have $W_M(V)\sbq\ker
(1\otimes\rho_{\Lambda^2,T^2}\otimes 1)$. It follows that\\
$1\otimes\rho_{\Lambda^2,T^2}\otimes 1:T(V)\otimes\Lambda^2(V)\otimes
T(V)\to T(V)$ induces a morphism of graded bimodules
$\rho_{M,T}:M(V)\to T(V)$, given by $[\xi\otimes x\wedge y\otimes\eta
]\mapsto\xi\otimes [x,y]\otimes\eta$ $\forall x,y\in V$ and
$\xi,\eta\in T(V)$. Moreover we have the exact sequence
$$M(V)\xrightarrow{\rho_{M,T}}T(V)\xrightarrow{\rho_{T,S}}S(V)\to 0.$$

Since $\rho_{M,T}$ is a morphism of graded bimodules, we may consider
its homogenous components, $\rho_{M^n,T^n}:M^n(V)\to T^n(V)$

\btm We have $W_M(V)=\ker (1\otimes\rho_{\Lambda^2,T^2}\otimes 1)$,
i.e. $\rho_{M,T}$ is injective and we have the exact sequence
$$0\to M(V)\xrightarrow{\rho_{M,T}}T(V)\xrightarrow{\rho_{T,S}}S(V)\to
0.$$
\etm
\pf We use induction on $n$ to prove that $\rho_{M^n,T^n}$ is
injective. If $n=0,1$ then $(T(V)\otimes\Lambda^2(V)\otimes
T(V))^n=0$ so $M^n(V)=0$, so there is nothing to prove. 

Before proving the induction step, we need some preliminary results.
\medskip

We have an action of the symmetric group $S_n$ on $T^n(V)$, given by
$$\tau (x_1\otimes\cdots\otimes
x_n)=x_{\tau^{-1}(1)}\otimes\cdots\otimes x_{\tau^{-1}(n)}.$$

For $1\leq i\leq n-1$ we denote by $\tau_i$ the transposition
$(i,i+1)\in S_n$ and we denote by $f_i:T^n(V)\to M^n(V)$ the linear
map given by $x_1\otimes\cdots\otimes x_n\mapsto
[x_1\otimes\cdots\otimes x_i\wedge x_{i+1}\otimes\cdots\otimes
x_n]$. (Here we repaced the $\otimes$ sign between $x_i$ and
$x_{i+1}$ by $\wedge$.)

\blm On $T^n(V)$ we have $\rho_{M^n,T^n}f_i=1-\tau_i$.
\elm

\pf We verify this relation on generators $\xi =x_1\otimes\cdots\otimes x_n$
of $T^n(V)$. We have
$$\begin{aligned}\rho_{M^n,T^n}f_i(\xi
)&=\rho_{M^n,T^n}([x_1\otimes\cdots\otimes x_i\wedge
x_{i+1}\otimes\cdots\otimes x_n])\\
{}&=x_1\otimes\cdots\otimes [x_i,x_{i+1}]\otimes\cdots\otimes x_n\\
{}&=x_1\otimes\cdots\otimes (x_i\otimes x_{i+1}-x_{i+1}\otimes
x_i)\otimes\cdots\otimes x_n\\
{}&=x_1\otimes\cdots\otimes x_n-x_1\otimes\cdots\otimes x_{i+1}\otimes
x_i\cdots\otimes x_n\\
{}&=x_1\otimes\cdots\otimes x_n-x_{\tau_i(1)}\otimes\cdots\otimes
x_{\tau_i(n)}=\xi -\tau_i(\xi ).\qquad\square\end{aligned}$$

\bco For every $1\leq i_1,\ldots,i_s\leq n-1$ in $T(V)$ we have
$$1-\tau_{i_s}\cdots\tau_{i_1}=\rho_{M^n,T^n}
\sum_{k=1}^sf_{i_k}\tau_{i_{k-1}}\cdots\tau_{i_1}.$$
\eco
\pf By Lemma 1.4, for $1\leq k\leq s$ in we have
$\rho_{M^n,T^n}f_{i_k}=1-\tau_{i_k}$ so
$\rho_{M^n,T^n}f_{i_k}\tau_{i_{k-1}}\cdots\tau_{i_1}
=(1-\tau_{i_k})\tau_{i_{k-1}}\cdots\tau_{i_1}$. Hence $\rho_{M^n,T^n}
\sum_{k=1}^sf_{i_k}\tau_{i_{k-1}}\cdots\tau_{i_1}$ is equal to the
telescoping sum
$\sum_{k=1}^s(1-\tau_{i_k})\tau_{i_{k-1}}\cdots\tau_{i_1}
=1-\tau_{i_s}\cdots\tau_{i_1}$. \qed

\blm (i) If $\tau\in S_n$ then there is a map $h_\tau :T(V)\to M(V)$
with $h_\tau=\sum_{k=1}^sf_{i_k}\tau_{i_{k-1}}\cdots\tau_{i_1}$
whenever $\tau =\tau_{i_s}\cdots\tau_{i_1}$. In particular, $h_1=0$
and $h_{\tau_i}=f_i$.

(ii) $h_{\sigma\tau}=h_\tau +h_\sigma\tau$ $\forall\sigma,\tau\in
S_n$.

(iii) On $T^n(V)$ we have $\rho_{M^k,T^k}h_\tau =1-\tau$
$\forall\tau\in S_n$.
\elm
\pf (i) We use the fact that $S_n$ is generated by
$\tau_1,\ldots,\tau_{n-1}$, with the relations $\tau_i^2=1$,
$\tau_i\tau_j=\tau_j\tau_i$ if $j-i\geq 2$ and
$(\tau_i\tau_{i+1})^3=1$. (See, e.g., [KT, Theorem 4.1, pag. 152].)

We consider the set of symbols
$A=\{\sigma_1,\ldots,\sigma_{n-1}\}$. Then $S_n$ is isomorphic to the
free monoid $(A\cup A^{-1})^*$ factored by the equivalence
relation $\sim$, generated by $\alpha\beta\sim\alpha\gamma\beta$ for
every $\alpha,\beta\in (A\cup A^{-1})^*$ and $\gamma$ of the form
$\gamma =\sigma_i\sigma_i^{-1}$ or $\sigma_i^{-1}\sigma_i$, $\gamma
=\sigma_i^2$, $\gamma=(\sigma_i\sigma_{i+1})^3$ or $\gamma
=\sigma_i\sigma_j\sigma_i^{-1}\sigma_j^{-1}$, with $j-i\geq 2$. For
any $\sigma\in (A\cup A^{-1})^*$ we denote by $[\sigma ]$ its class in
$(A\cup A^{-1})^*_{/\sim}$. If $\psi :(A\cup A^{-1})^*\to S_n$ is the
morphism of monoids given by $\sigma_i\mapsto\tau_i$ and
$\sigma_i^{-1}\mapsto\tau_i^{-1}=\tau_i$ then $\psi$ induces an
isomorphism $\tilde\psi :(A\cup A^{-1})^*_{/\sim}\to S_n$, given by
$\tilde\psi ([\sigma ])=\psi (\sigma )$ $\forall\sigma\in (A\cup
A^{-1})^*$.

For every $\sigma =\sigma_s^{\pm 1}\cdots\sigma_1^{\pm 1}\in (A\cup
A^{-1})^*$ define the map $g_\sigma
=\sum_{k=1}^sf_{i_k}\tau_{i_{k-1}}\cdots\tau_{i_1}$. (If $\sigma =1$
then $s=0$ so $g_1:=0$.)

Note that if $\alpha=\sigma_s^{\pm 1}\cdots\sigma_{t+1}^{\pm 1}$ and
$\beta =\sigma_t^{\pm 1}\cdots\sigma_1^{\pm 1}$ then $\psi (\beta
)=\tau_t\cdots\tau_1$ so

$$\begin{aligned}g_{\alpha\beta}&
=\sum_{k=1}^sf_{i_k}\tau_{i_{k-1}}\cdots\tau_{i_1}
=\sum_{k=1}^tf_{i_k}\tau_{i_{k-1}}\cdots\tau_{i_1}
+\left(\sum_{k=t+1}^sf_{i_k}\tau_{i_{k-1}}\cdots\tau_{i_{t+1}}\right)
\tau_{i_t}\cdots\tau_{i_1}\\
&=g_\beta+g_\alpha\psi (\beta ).\end{aligned}$$

We now prove that if $\sigma\sim\sigma'$ then $g_\sigma
=g_{\sigma'}$. It suffices to take the case when
$\sigma=\alpha\gamma\beta$ and $\sigma'=\alpha\beta$, with
$\alpha,\beta\in (A\cup A^{-1})^*$ and $\gamma$ is of the form
$\sigma_i\sigma_i^{-1}$, $\sigma_i^{-1}\sigma_i$, $\sigma_i^2$,
$(\sigma_i\sigma_{i+1})^3$ or
$\sigma_i\sigma_j\sigma_i^{-1}\sigma_j^{-1}$, with $j-i\geq 2$. Note
that in all these cases we have $\psi (\gamma )=1$. (We have
$\tau_i^2=1$, $(\tau_i\tau_{i+1})^3=1$ and, if $j-i\geq 2$, then
$\tau_i\tau_j\tau_i\tau_j=\tau_i^2\tau_j^2=1$.)

The relation $g_{\alpha\beta}=g_{\alpha\gamma\beta}$ writes as
$g_\beta +g_\alpha\psi (\beta )=g_\beta +g_{\alpha\gamma}\psi (\beta
)$ so it suffices to prove that $g_\alpha=g_{\alpha\gamma}$. But $\psi
(\gamma )=1$ so $g_{\alpha\gamma}=g_\gamma +g_\alpha\psi (\gamma
)=g_\gamma +g_\alpha$. Hence we must prove that $g_\gamma =0$. We
prove that $g_\gamma (\eta )=0$ for $\eta =x_1\otimes\cdots\otimes
x_n$.

If $\gamma =\sigma_i^{\pm 1}\sigma_i^{\pm 1}$, which includes the
cases $\gamma =\sigma_i\sigma_i^{-1}$, $\sigma_i^{-1}\sigma_i$ and
$\sigma_i^2$, we have

$$\begin{aligned}g_\gamma (\eta )&=f_i(\eta )+f_i(\tau_i(\eta
))=f_i(x_1\otimes\cdots\otimes x_n)+f_i(x_1\otimes\cdots\otimes
x_{i+1}\otimes x_i\otimes\cdots\otimes x_n)\\
&=x_1\otimes\cdots\otimes x_i\wedge x_{i+1}\otimes\cdots\otimes
x_n+x_1\otimes\cdots\otimes x_{i+1}\wedge x_i\otimes\cdots\otimes
x_n=0.\end{aligned}$$

Let now $\gamma
=(\tau_i\tau_{i+1})^3=\tau_i\tau_{i+1}\tau_i\tau_{i+1}\tau_i\tau_{i+1}$. We
have $\eta=\eta'\otimes x\otimes y\otimes z\otimes\eta''$, where
$\eta'=x_1\otimes\cdots\otimes x_{i-1}$,
$\eta''=x_{i+3}\otimes\cdots\otimes x_n$ and
$(x,y,z)=(x_i,x_{i+1},x_{i+2})$. Note that when we apply succesively
the transpositions $\tau_i=(i,i+1)$ and $\tau_{i+1}=(i+1,i+2)$ to
$\eta$ only the factors $x,y,z$ of $\eta$ are permuted, while the
factors $\eta'$ and $\eta''$ are unchanged. The factors on the
positions $i$, $i+1$ and $i+2$ in $\eta$ are $x,y,z$; in
$\tau_{i+1}(\eta )$ they are $x,z,y$; in $\tau_i\tau_{i+1}(\eta )$
they are $z,x,y$; in $\tau_{i+1}\tau_i\tau_{i+1}(\eta )$ they are
$z,y,x$; in $\tau_i\tau_{i+1}\tau_i\tau_{i+1}(\eta )$ they are
$y,z,x$; and in $\tau_{i+1}\tau_i\tau_{i+1}\tau_i\tau_{i+1}(\eta )$
they are $y,x,z$. Therefore

$$\begin{aligned}g_\gamma (\eta )&=f_{i+1}(\eta )+f_i\tau_{i+1}(\eta
)+f_{i+1}\tau_i\tau_{i+1}(\eta )+f_i\tau_{i+1}\tau_i\tau_{i+1}(\eta )\\
&\qquad\qquad\qquad
+f_{i+1}\tau_i\tau_{i+1}\tau_i\tau_{i+1}(\eta
)+f_i\tau_{i+1}\tau_i\tau_{i+1}\tau_i\tau_{i+1}(\eta )\\
&=f_{i+1}(\eta'\otimes x\otimes y\otimes
z\otimes\eta'')+f_i(\eta'\otimes x\otimes z\otimes
y\otimes\eta'')\\
&\qquad\qquad\qquad +f_{i+1}(\eta'\otimes z\otimes
x\otimes y\otimes\eta'')+f_i(\eta'\otimes z\otimes y\otimes
x\otimes\eta'')\\
&\qquad\qquad\qquad +f_{i+1}(\eta'\otimes y\otimes z\otimes
x\otimes\eta'')+f_i(\eta'\otimes y\otimes x\otimes z\otimes\eta'')\\
&= [\eta'\otimes x\otimes y\wedge z\otimes\eta'']+[\eta'\otimes
x\wedge z\otimes y\otimes\eta'']+[\eta'\otimes z\otimes x\wedge
y\otimes\eta'']\\
&\qquad +[\eta'\otimes z\wedge y\otimes
x\otimes\eta'']+[\eta'\otimes y\otimes z\wedge
x\otimes\eta'']+[\eta'\otimes y\wedge x\otimes
z\otimes\eta''].\end{aligned}$$
Thus $g_\gamma (\eta )=[\eta'\otimes\xi\otimes\eta'']$, where
$$\begin{aligned}\xi&=x\otimes y\wedge z+x\wedge z\otimes y+z\otimes
x\wedge y+z\wedge y\otimes x+y\otimes z\wedge x+y\wedge x\otimes z\\ 
& =x\otimes y\wedge z-z\wedge x\otimes y+z\otimes x\wedge y-y\wedge
z\otimes x+y\otimes z\wedge x-x\wedge y\otimes z\\
&=[x,y\wedge z]+[y,z\wedge x]+[z,x\wedge y].\end{aligned}$$
We have $\xi\in W_M(V)$ so $\eta'\otimes\xi\otimes\eta''\in W_M(V)$ and so
$g_\gamma (\eta )=[\eta'\otimes\xi\otimes\eta'']=0$.

Let now $\gamma =\sigma_i\sigma_j\sigma_i\sigma_j$. We have $\eta
=\eta'\otimes x\otimes y\otimes\xi\otimes z\otimes t\otimes\eta''$,
where $\eta'=x_1\otimes\cdots\otimes x_{i-1}$,
$\xi=x_{i+2}\otimes\cdots\otimes x_{j-1}$,
$\eta''=x_{j+2}\otimes\cdots\otimes x_n$, $(x,y)=(x_i,x_{i+1})$ and
$(z,t)=(x_j,x_{j+1})$. Note that $\tau_i$ permutes the factors $x$
and $y$ of $\eta$ and leaves all the other factors unchanged, while
$\tau_j$ permutes $z$ and $t$ and leaves all the other factors
unchanged. We get $\tau_j(\eta )=\eta'\otimes x\otimes
y\otimes\xi\otimes t\otimes z\otimes\eta''$, $\tau_i\tau_j(\eta
)=\eta'\otimes y\otimes x\otimes\xi\otimes t\otimes z\otimes\eta''$
and $\tau_j\tau_i\tau_j(\eta )=\eta'\otimes y\otimes
x\otimes\xi\otimes z\otimes t\otimes\eta''$. Then

$$\begin{aligned}g_\gamma (\eta )&=f_j(\eta )+f_i\tau_j(\eta
)+f_j\tau_i\tau_j(\eta )+f_i\tau_j\tau_i\tau_j(\eta )\\
&=f_j(\eta'\otimes x\otimes y\otimes\xi\otimes z\otimes
t\otimes\eta'')+f_i(\eta'\otimes x\otimes y\otimes\xi\otimes
t\otimes z\otimes\eta'')\\
&\qquad\qquad +f_j(\eta'\otimes y\otimes x\otimes\xi\otimes t\otimes
z\otimes\eta'')+f_i(\eta'\otimes y\otimes x\otimes\xi\otimes
z\otimes t\otimes\eta'')\\
&=[\eta'\otimes x\otimes
y\otimes\xi\otimes z\wedge t\otimes\eta'']+[\eta'\otimes x\wedge
y\otimes\xi\otimes t\otimes z\otimes\eta'']\\
&\qquad\qquad +[\eta'\otimes y\otimes x\otimes\xi\otimes t\wedge
z\otimes\eta'']+[\eta'\otimes y\wedge x\otimes\xi\otimes z\otimes
t\otimes\eta''].\end{aligned}$$
Thus $g_\gamma (\eta )=[\eta'\otimes\xi'\otimes\eta'']$, where
$$\begin{aligned}\xi'&=x\otimes y\otimes\xi\otimes z\wedge t+x\wedge
y\otimes\xi\otimes t\otimes z+y\otimes x\otimes\xi\otimes t\wedge
z+y\wedge x\otimes\xi\otimes z\otimes t\\
&=x\otimes y\otimes\xi\otimes z\wedge t+x\wedge
y\otimes\xi\otimes t\otimes z-y\otimes x\otimes\xi\otimes z\wedge
t-x\wedge y\otimes\xi\otimes z\otimes t\\
&=[x,y]\otimes\xi\otimes z\wedge t-x\wedge
y\otimes\xi\otimes [z,t].\end{aligned}$$
We have $\xi'\in W_M(V)$ so $\eta'\otimes\xi'\otimes\eta''\in
W_M(V)$ and $g_\gamma (\eta )=[\eta'\otimes\xi'\otimes\eta'']=0$.

Since the map $\sigma\mapsto g_\sigma$, defined on $(A\cup A^{-1})^*$,
is invariant to the equivalence relation $\sim$, it induces a map
defined on $(A\cup A^{-1})^*_{/\sim}$, given by $[\sigma ]\mapsto
g_\sigma$. Since $\bar\psi :(A\cup A^{-1})^*_{/\sim}\to S_n$ is an
isomorphism we get a map $\tau\mapsto h_\tau$, where $h_\tau=g_\sigma$
for any $\sigma\in (A\cup A^{-1})^*$ such that $\tau =\bar\psi
([\sigma ])=\psi (\sigma )$. If $\tau =\tau_{i_s}\cdots\tau_{i_1}$
then $\tau =\psi (\sigma )$, with $\sigma
=\sigma_{i_s}\cdots\sigma_{i_1}$. Hence $h_\tau=g_\sigma
=\sum_{k=1}^sf_k\tau_{k-1}\cdots\tau_1$, as claimed.

(ii) Let $\alpha,\beta\in (A\cup A^{-1})^*$ with $\sigma =\psi (\alpha
)$ and $\tau =\psi (\beta )$ so that $\sigma\tau =\psi (\alpha\beta
)$. Then $h_\sigma =g_\alpha$, $h_\tau =g_\beta$ and
$h_{\sigma\tau}=g_{\alpha\beta}=g_\beta +g_\alpha\psi (\beta )=h_\tau
+h_\sigma\tau$.

(iii) If we write $\tau=\tau_{i_s}\cdots\tau_1$ then our result is
just Corollary 1.5.\qed

{\bf Proof of the induction step.} We must prove that if $[\eta
]\in\ker\rho_{M^n,T^n}$ then $[\eta ]=0$. Note that $\eta$ is a finite
linear combination of products of the form\\
$v_{i_1}\otimes\cdots\otimes v_{i_k}\wedge
v_{i_{k+1}}\otimes\cdots\otimes v_{i_n}$, with $i_1,\ldots,i_n\in I$
and $i_k<i_{k+1}$. We denote by $J=\{ j_1,\ldots,j_m\}$ with
$j_1,\ldots,j_m\in I$, $j_1<\cdots <j_m$, the set of all indices $i\in
I$ such that $v_i$ is one of the factors $v_{i_h}$ from one of the
products in the linear combination that gives $\eta$. We will prove
our result by induction on $m$. If $m=0$ there is nothing to be
proved. Suppose that $m\geq 1$. Let $J'=\{ j_1,\ldots,j_{m-1}\}$.

We have $\eta\in W$, where $W\sbq (T(V)\otimes\Lambda^2(V)\otimes
T(V))^n$ is spanned by $v_{i_1}\otimes\cdots\otimes v_{i_k}\wedge
v_{i_{k+1}}\otimes\cdots\otimes v_{i_n}$, with $(i_1,\ldots,i_n,;k)\in
A$, where $A$ is the set of all $(i_1,\ldots,i_n,;k)$ with
$i_1,\ldots,i_n\in J$, $1\leq k\leq n-1$ and $i_k<i_{k+1}$. We also
denote by $U\sbq T^n(V)$ the space generated by
$v_{i_1}\otimes\cdots\otimes v_{i_n}$, with $(i_1,\ldots,i_n)\in
B:=J^n$.

We have $A=A_1\sqcup A_2$ where $A_1$ is the set of all
$(i_1,\ldots,i_n;k)\in A$ with $i_1,\ldots,i_n\in J'$ and $A_2$ is the
set of those where at least one of $i_1,\ldots,i_n$ is
$j_m$. Similarly, $B=B_1\sqcup B_2$ where $B_1=J'^n$ and
$B_2=J^n\setminus J'^n$, i.e. $B_2$ is the set of all
$(i_1,\ldots,i_n)\in J^n$ such that at least one of $i_1,\ldots,i_n$
is $j_m$. For $\alpha =1,2$ we denote by $W_\alpha$ the subspace of
$W$ spanned by $v_{i_1}\otimes\cdots\otimes v_{i_k}\wedge
v_{i_{k+1}}\otimes\cdots\otimes v_{i_n}$ with $(i_1,\ldots,i_n,;k)\in
A_\alpha$ and by $U_\alpha$ the subspace of $U$ spanned by
$v_{i_1}\otimes\cdots\otimes v_{i_n}$ with $(i_1,\ldots,i_n)\in
B_\alpha$. Then from $A=A_1\sqcup A_2$
and $B=B_1\sqcup B_2$ we deduce that $W=W_1\oplus
W_2$ and $U=U_1\oplus U_2$.

We have $(1\otimes\rho_{\Lambda^2,T^2}\otimes
1)(v_{i_1}\otimes\cdots\otimes v_{i_k}\wedge
v_{i_{k+1}}\otimes\cdots\otimes v_{i_n})=v_{i_1}\otimes\cdots\otimes
v_{i_n}-v_{i_1}\otimes\cdots\otimes v_{i_{k+1}}\otimes
v_{i_k}\otimes\cdots\otimes v_{i_n}$. If $(i_1,\ldots,i_n;k)\in A$,
$A_1$ or $A_2$ then both $(i_1,\ldots,i_n)$ and
$(i_1,\ldots,i_{k+1},i_k,\ldots,i_n)$ belong to $B$, $B_1$ or $B_2$
and so $(1\otimes\rho_{\Lambda^2,T^2}\otimes
1)(v_{i_1}\otimes\cdots\otimes v_{i_k}\wedge
v_{i_{k+1}}\otimes\cdots\otimes v_{i_n})\in U$, $U_1$ or $U_2$,
respectively. It follows that $(1\otimes\rho_{\Lambda^2,T^2}\otimes
1)(W)\sbq U$ and $(1\otimes\rho_{\Lambda^2,T^2}\otimes 1)(W_\alpha
)\sbq U_\alpha$ for $\alpha =1,2$. Equivalently, if $\xi\in W$, $W_1$
or $W_2$ then $\rho_{M^n,T^n}([\xi
])=(1\otimes\rho_{\Lambda^2,T^2}\otimes 1)(\xi )\in U$, $U_1$ or
$U_2$, respectively.

We define $\psi :U_2\to M^n(V)$ on elements on the basis as
follows. If $\xi =v_{i_1}\otimes\cdots\otimes v_{i_n}$ with
$(i_1,\ldots,i_n)\in B_2$ and $l$ is the smallest index vith $i_l=j_m$
then $\psi (\xi ):=h_{\sigma_l}(\xi )$, where $\sigma_l\in S_n$ is the
cyclic permutation $(1,2,\ldots,l)$.

After these preliminaries, we start our proof of the induction step.

Since $W=W_1\oplus W_2$ we have $\eta =\eta_1+\eta_2$, with
$\eta_\alpha\in W_\alpha$. Then $\rho_{M^n,T^n}([\eta_\alpha ])\in
U_\alpha$. Since
$\rho_{M^n,T^n}([\eta_1])+\rho_{M^n,T^n}([\eta_2])=\rho_{M^n,T^n}([\eta
])=0$ and the sum $U_1+U_2$ is direct, this implies that
$\rho_{M^n,T^n}([\eta_1])=\rho_{M^n,T^n}([\eta_2])=0$. Since
$\eta_1\in U_1$, it can be written in terms of only
$v_i$ with $i\in J'$. Since $|J'|=m-1$, by the induction hypothesis,
$\rho_{M^n,T^n}([\eta_1])=0$ implies $[\eta_1]=0$ so $[\eta
]=[\eta_2]$. So we have reduced to the case when $\eta\in U_2$. Then
$\eta$ writes as
$$\eta =\sum a_{i_1,\ldots,i_n;k}v_{i_1}\otimes\cdots\otimes v_{i_k}\wedge
v_{i_{k+1}}\otimes\cdots\otimes v_{i_n},$$
where the sum is take over $(i_1,\ldots,i_n;k)\in A_2$ and
$a_{i_1,\ldots,i_n;k}\in F$. Since $[v_{i_1}\otimes\cdots\otimes v_{i_k}\wedge
v_{i_{k+1}}\otimes\cdots\otimes
v_{i_n}]=f_k(v_{i_1}\otimes\cdots\otimes v_{i_n})$ we have
$$[\eta ]=\sum a_{i_1,\ldots,i_n;k}f_k(v_{i_1}\otimes\cdots\otimes
v_{i_n}).$$

By Lemma 1.4, on $T^n(V)$ we have $\rho_{M^n,T^n}f_k=1-\tau_k$. It
follows that
$$0=\rho_{M^n,T^n}[\eta ]=\sum
a_{i_1,\ldots,i_n;k}(v_{i_1}\otimes\cdots\otimes
v_{i_n}-\tau_k(v_{i_1}\otimes\cdots\otimes v_{i_n})).$$
But for every $(i_1,\ldots,i_n;k)\in A_2$ we have $(i_1,\ldots,i_n)\in
B_2$, which implies that also $(i_{\tau_k(1)},\ldots,i_{\tau_k(n)})\in
B_2$. Thus both $v_{i_1}\otimes\cdots\otimes v_{i_n}$ and
$\tau_k(v_{i_1}\otimes\cdots\otimes v_{i_n})$ belong to $U_2$ so we
can apply $\psi$ to the formula above. We get
$$0=\sum a_{i_1,\ldots,i_n;k}(\psi (v_{i_1}\otimes\cdots\otimes
v_{i_n})-\psi\tau_k(v_{i_1}\otimes\cdots\otimes v_{i_n})).$$

Let $(i_1,\ldots,i_n;k)\in A_2$. We have $\psi
(v_{i_1}\otimes\cdots\otimes
v_{i_n})=h_{\sigma_l}(v_{i_1}\otimes\cdots\otimes v_{i_n})$ and
$\psi\tau_k(v_{i_1}\otimes\cdots\otimes
v_{i_n})=h_{\sigma_{l'}}\tau_k(v_{i_1}\otimes\cdots\otimes v_{i_n})$,
where $l$ is the smallest index with $i_l=j_m$ and $l'$ is the
smallest index with $i_{\tau_k(l')}=j_m$. Since $i_k<i_{k+1}\leq j_m$,
we cannot have $l=k$. Since $(i_{\tau_k(1)},\ldots,i_{\tau_k(n)})=
(i_1,\ldots,i_{k-1},i_{k+1},i_k,i_{k+2},\ldots,i_n)$, if $l\leq k-1$
or $l\geq k+2$ then $l'=l$. If $l=k+1$ then $l'=k$.

Let $\tau =\sigma_{l'}\tau_k\sigma_l^{-1}$ so that
$\tau\sigma_l=\sigma_{l'}\tau_k$. Then, by Lemma 1.6(ii), the relation
$h_{\tau\sigma_l}=h_{\sigma_{l'}\tau_k}$ writes as
$$h_{\sigma_l}+h_\tau\sigma_l =h_{\tau_k}+h_{\sigma_{l'}}\tau_k
=f_k+h_{\sigma_{l'}}\tau_k$$
so $f_k=h_\tau\sigma_l+(h_{\sigma_l}-h_{\sigma_{l'}}\tau_k)$. Since
$h_{\sigma_l}(v_{i_1}\otimes\cdots\otimes v_{i_n})=\psi
(v_{i_1}\otimes\cdots\otimes v_{i_n})$ and
$h_{\sigma_{l'}}\tau_k(v_{i_1}\otimes\cdots\otimes
v_{i_n})=\psi\tau_k(v_{i_1}\otimes\cdots\otimes v_{i_n})$, this
implies that
$$f_k(v_{i_1}\otimes\cdots\otimes
v_{i_n})=[\zeta_{i_1,\ldots,i_n;k}]+(\psi (v_{i_1}\otimes\cdots\otimes
v_{i_n})-\psi\tau_k(v_{i_1}\otimes\cdots\otimes v_{i_n})),$$
where $[\zeta_{i_1,\ldots,i_n;k}]
=h_\tau\sigma_l(v_{i_1}\otimes\cdots\otimes v_{i_n})$.

It follows that
$$\begin{aligned} \left[\eta \right] &=\sum
a_{i_1,\ldots,i_n;k}f_k(v_{i_1}\otimes\cdots\otimes v_{i_n})\\
{}&=\sum a_{i_1,\ldots,i_n;k}[\zeta_{i_1,\ldots,i_n;k}]+\sum
a_{i_1,\ldots,i_n;k}(\psi (v_{i_1}\otimes\cdots\otimes
v_{i_n})-\psi\tau_k(v_{i_1}\otimes\cdots\otimes v_{i_n}))\\
{}&=\sum a_{i_1,\ldots,i_n;k}[\zeta_{i_1,\ldots,i_n;k}].\end{aligned}$$

We now claim that if $(i_1,\ldots,i_n;k)\in A_2$ then
$[\zeta_{i_1,\ldots,i_n;k}]=v_{j_m}\otimes [\zeta'_{i_1,\ldots,i_n;k}]$
for some $\zeta'_{i_1,\ldots,i_n;k}\in (T(V)\otimes\Lambda^2(V)\otimes
T(V))^{n-1}$.

If $l=k+1$ then $l'=k$ so $\tau =\sigma_k\tau_k\sigma_{k+1}^{-1}$. But
in $S_n$ we have $(1,2,\ldots,k)(k,k+1)=(1,2,\ldots,k+1)$,
i.e. $\sigma_k\tau_k=\sigma_{k+1}$, so $\tau =1$, which implies
$h_\tau =0$, so $[\zeta_{i_1,\ldots,i_n;k}]
=h_\tau\sigma_l(v_{i_1}\otimes\cdots\otimes v_{i_n})=0$ and we may
take $\zeta'_{i_1,\ldots,i_n;k}=0$.

Suppose now that $l\leq k-1$ or $l\geq k+2$, so that $l'=l$. We have
$\sigma_l(v_{i_1}\otimes\cdots\otimes
v_{i_n})=v_{i'_1}\otimes\cdots\otimes v_{i'_n}$, with
$i'_h=i_{\sigma_l^{-1}(h)}$. But $\sigma_l(l)=1$ so
$i'_1=i_{\sigma_l^{-1}(1)}=i_l=j_m$. So
$\sigma_l(v_{i_1}\otimes\cdots\otimes v_{i_n})=v_{j_m}\otimes
v_{i'_2}\otimes\cdots\otimes v_{i'_n}$.

Since $l'=l$ we have $\tau =\sigma_l\tau_k\sigma_l^{-1}
=\sigma_l(k,k+1)\sigma_l^{-1} =(\sigma_l(k),\sigma_l(k+1))$.

If $l\leq k-1$ then $\sigma_l(k)=k$ and $\sigma_l(k+1)=k+1$ so $\tau
=(k,k+1)=\tau_k$, so $h_\tau =f_k$. Thus
$[\zeta_{i_1,\ldots,i_n;k}]=h_\tau\sigma_l(v_{i_1}\otimes\cdots\otimes
v_{i_n})=f_k(v_{j_m}\otimes v_{i'_1}\otimes\cdots\otimes
v_{i'_n})=[v_{j_m}\otimes\zeta'_{i_1,\ldots,i_n;k}]$, where
$\zeta'_{i_1,\ldots,i_n;k}=v_{i'_2}\otimes\cdots\otimes v_{i'_k}\wedge
v_{i'_{k+1}}\otimes\cdots\otimes v_{i'_n}$. (Note that $l\leq k-1$
implies $k\geq 2$.)

If $l\geq k+2$ then $\sigma_l(k)=k+1$ and $\sigma_l(k+1)=k+2$ so $\tau
=(k+1,k+2)=\tau_{k+1}$, so $h_\tau =f_{k+1}$. Then, by the same
reasoning from the case $l\leq k-1$, we get
$[\zeta_{i_1,\ldots,i_n;k}]=[v_{j_m}\otimes\zeta'_{i_1,\ldots,i_n;k}]$, where
$\zeta'_{i_1,\ldots,i_n;k}=v_{i'_2}\otimes\cdots\otimes
v_{i'_{k+1}}\wedge v_{i'_{k+2}}\otimes\cdots\otimes v_{i'_n}$.

Since $[\zeta_{i_1,\ldots,i_n;k}]=v_{j_m}\otimes
[\zeta'_{i_1,\ldots,i_n;k}]$ we have $[\eta ]=\sum
a_{i_1,\ldots,i_n;k}[\zeta_{i_1,\ldots,i_n;k}]=v_{j_m}\otimes
[\eta']$, where $\eta'=\sum
a_{i_1,\ldots,i_n;k}\zeta'_{i_1,\ldots,i_n;k}$. Then
$0=\rho_{M^n,T^n}([\eta ])=\rho_{M^n,T^n}(v_{j_m}\otimes [\eta'])=
v_{j_m}\otimes\rho_{M^{n-1},T^{n-1}}([\eta'])$. It follows that
$\rho_{M^{n-1},T^{n-1}}([\eta'])=0$. But, by the induction hypothesis,
$\rho_{M^{n-1},T^{n-1}}$ is injective so we have $[\eta']=0$. It
follows that $[\eta ]=v_{j_m}\otimes [\eta']=0$. \qed

{\bf Remark} When $n\geq 2$ we have $(T(V)\otimes\Lambda^2(V)\otimes
T(V))^2=\Lambda^2(V)$ and $W_M^2(V)=0$. (All generators of $W_M(V)$
have degrees $\geq 3$.) Thus $M^2(V)=\Lambda^2(V)$ and $\rho_{M^2,T^2}$
coincides with $\rho_{\Lambda^2,T^2}$. Hence $0\to
M^2(V)\xrightarrow{\rho_{M^2,T^2}} T^2(V)\xrightarrow{\rho_{T^2,S^2}}
S^2(V)\to 0$ coincides with the short exact sequence from Proposition
1.1. 

\section{The algebra $S'(V)$}

\bdf We define the algebra $S'(V)$ as $S'(V)=T(V)/W_{S'}(V)$, where
$W_{S'}(V)$ is the bilateral ideal of $T(V)$ generated by $x\otimes
y\otimes z-y\otimes z\otimes x$, with $x,y,z\in V$.

If $x_1\ldots,x_n\in V$ then we denote by $x_1\odot\cdots\odot x_n$
the class of $x_1\otimes\cdots\otimes x_n$ in $S'(V)$. So
$(S'(V),+,\odot )$ is an algebra.
\edf

Since $W_{S'}(V)$ is a homogeneous ideal, $S'(V)$ is a graded
algebra. For $n\geq 0$ we denote by $S'^n(V)$ the homogeneous
component of degree $n$ of $S'(V)$. We have
$S'^n(V)=T^n(V)/W^n_{S'}(V)$, where $W^n_{S'}(V)=W_{S'}(V)\cap
T^n(V)$. 

Since the generators of $W_{S'}(V)$ have degree $3$, for $n\leq 2$ we
have $W^n_{S'}(V)=0$ so $S'^n(V)=T^n(V)$. 

Note that $\rho_{T,S}(x\otimes y\otimes z-y\otimes z\otimes
x)=xyz-yzx=0$ so $x\otimes y\otimes z-y\otimes z\otimes
x\in\ker\rho_{T,S}$ $\forall x,y,z\in V$. It follows that
$W_{S'}(V)\sbq\ker\rho_{S,T}$. Therefore $\rho_{T,S}$ induces a
surjective morphism of algebras defined on
$T(V)/W_{S'}(V)=S'(V)$. Namely, we have:

\bpr There is a surjective morphism of algebras $\rho_{S',S}:S'(V)\to
S(V)$ given by $x_1\odot\cdots\odot x_n\mapsto x_1\cdots x_n$.
\epr

\bpr The subspace $W_{S'}^n(V)$ of $T(V)$ is spanned by
$x_1\otimes\cdots\otimes x_n-x_{\sigma (1)}\otimes\cdots\otimes
x_{\sigma (n)}$, with $x_1,\ldots,x_n\in V$ and $\sigma\in A_n$.
\epr
\pf The bilateral ideal $W_{S'}(V)$ of $T(V)$ is generated by
$f(x,y,z)$, where $f:V^3\to T^3(V)$ is given by $(x,y,z)\mapsto
x\otimes y\otimes z-y\otimes z\otimes x$. Therefore $W_{S'}^n(V)$ is
spanned by
\begin{multline*}x_1\otimes\cdots\otimes x_{i-1}\otimes
f(x_i,x_{i+1},x_{i+2})\otimes x_{i+3}\otimes\cdots\otimes x_n\\
=x_1\otimes\cdots\otimes x_n-x_1\otimes\cdots\otimes
x_{i+1}\otimes x_{i+2}\otimes x_i\otimes\cdots\otimes x_n\\
=x_1\otimes\cdots\otimes x_n-x_{\sigma_i(1)}\otimes\cdots\otimes
x_{\sigma_i(n)},\end{multline*}
with $x_1,\ldots x_n\in V$ and $1\leq i\leq n-2$, where $\sigma_i\in
S_n$ is the cycle $(i,i+1,i+2)$.

Hence in $S'^n(V)$ we have $x_1\odot\cdots\odot
x_n=x_{\sigma_i(1)}\odot\cdots\odot x_{\sigma_i(n)}$ $\forall
x_1,\ldots,x_n\in V$ and $1\leq i\leq n-2$. But the cycles $\sigma_i$
generate the the alternating group $A_n$ so in $S'^n(V)$ we have
$x_1\odot\cdots\odot x_n=x_{\sigma (1)}\odot\cdots\odot x_{\sigma
(n)}$ $\forall\sigma\in A_n$. Hence $W_{S'}^n(V)$ contains
$x_1\otimes\cdots\otimes x_n-x_{\sigma (1)}\otimes\cdots\otimes
x_{\sigma (n)}$ $\forall x_1,\ldots,x_n\in V$ and $\sigma\in A_n$,
which are more general than the original generators, where $\sigma
=\sigma_i$ for some $1\leq i\leq n-2$. \qed

\bpr (i) If $n\geq 2$ then we have a linear map $c:S'^n(V)\to S'^n(V)$
given by $x_1\odot\cdots\odot x_n\mapsto x_2\odot x_1\odot
x_3\odot\cdots\odot x_n$, $\forall x_1,\ldots,x_n\in V$.

(ii) We have $c^2=1$ and $\rho_{S',S}c=\rho_{S',S}$

(iii) If $x_1,\ldots,x_n\in V$ and $\sigma\in S_n$ then
$$x_{\sigma (1)}\odot\cdots\odot x_{\sigma
(n)}=\begin{cases}x_1\odot\cdots\odot x_n&\text{ if }\sigma\in
A_n\\ c(x_1\odot\cdots\odot x_n)&\text{ if }\sigma\in S_n\setminus
A_n\end{cases}.$$

(iv) If there are $i<j$ with $x_i=x_j$ then $c(x_1\odot\cdots\odot
x_n)=x_1\odot\cdots\odot x_n$. Consequently, $x_{\sigma
(1)}\odot\cdots\odot x_{\sigma (n)}=x_1\odot\cdots\odot x_n$ holds
regardless of the parity of $\sigma$.
\epr
\pf (i) We define $\bar c:T^n(V)\to T^n(V)$ by
$x_1\otimes\cdots\otimes x_n\mapsto x_2\otimes x_1\otimes
x_3\otimes\cdots\otimes x_n=x_{\tau (1)}\otimes\cdots\otimes x_{\tau
(n)}$, where $\tau\in S_n$ is the transposition $(1,2)$. To prove that
$\bar c$ induces the morphism $c:S'^n(V)\to S'^n(V)$ given by
$x_1\odot\cdots\odot x_n\mapsto x_2\odot x_1\odot x_3\odot\cdots\odot
x_n$, one must prove that $\bar c(W_{S'}^n(V))\sbq W_{S'}^n(V)$.

Let $\xi =x_1\otimes\cdots\otimes x_n-x_{\sigma
(1)}\otimes\cdots\otimes x_{\sigma (n)}$ be a generator of $W_{S'}^n$,
with $x_1,\ldots,x_n\in V$ and $\sigma\in A_n$. Then $\bar c(\xi
)=x_{\tau (1)}\otimes\cdots\otimes x_{\tau (n)}-x_{\sigma\tau
(1)}\otimes\cdots\otimes x_{\sigma\tau (n)}$. If we denote
$y_i=x_{\tau (i)}$, so that $x_i=y_{\tau^{-1}(i)}=y_{\tau (i)}$, then
$x_{\sigma\tau (i)}=y_{\tau\sigma\tau (i)}$. Hence $\bar c(\xi
)=y_1\otimes\cdots\otimes y_n-y_{\sigma'(1)}\otimes\cdots\otimes
y_{\sigma'(n)}$, where $\sigma'=\tau\sigma\tau$. But $\sigma\in A_n$,
which implies that $\sigma'\in A_n$ and so $\bar c(\xi )\in
W_{S'}^n(V)$.

(ii) Let $\xi =x_1\odot\cdots\odot x_n$. Applying twice $c$ to $\xi$
permutes the first two factors of $\xi$ twice so we have $c^2(\xi
)=\xi$. We have $\rho_{S',S}c(\xi )=x_2x_1x_3\cdots x_n=x_1\cdots
x_n=\rho_{S',S}(\xi )$.

(iii) We have $c(x_1\odot\cdots\odot x_n)=y_1\odot\cdots\odot y_n$,
where $y_i=x_{\tau (i)}$, with $\tau =(1,2)\in S_n$. Then $x_i=y_{\tau
(i)}$.

If $\sigma\in A_n$ then $x_{\sigma (1)}\odot\cdots\odot x_{\sigma
(n)}=x_1\odot\cdots\odot x_n$ follows from Proposition 2.2. If
$\sigma\notin A_n$ then note that $x_{\sigma (i)}=y_{\tau\sigma (i)}$
and, since $\tau,\sigma\notin A_n$, we have $\tau\sigma\in
A_n$. Therefore $x_{\sigma (1)}\odot\cdots\odot x_{\sigma
(n)}=y_{\tau\sigma (1)}\odot\cdots\odot y_{\tau\sigma
(n)}=y_1\odot\cdots\odot y_n=c(x_1\odot\cdots\odot x_n)$.

(iv) Let $\tau\in S_n$, $\tau =(i,j)$. Since $x_i=x_j$, permuting the
factors $x_i$ and $x_j$ has no effect on te product
$x_1\odot\cdots\odot x_n$. Hence $x_{\tau (1)}\odot\cdots\odot x_{\tau
(n)}=x_1\odot\cdots\odot x_n$. But $\tau\in S_n\setminus A_n$ so, by
(iii), $x_{\tau (1)}\odot\cdots\odot x_{\tau
(n)}=c(x_1\odot\cdots\odot x_n)$. Hence the conclusion. \qed

We now produce a basis for $S'^n(V)$. For this purpose we need the
following elementary result.

\blm Let $U$ be a vector space with the basis $(u_\alpha )_{\alpha\in
A}$. Let $\sim$ be an equivalence relation on $A$ and let $B$ be a
set of representatives for $A_{/\sim}$.

Let $W\sbq U$ be the subspace generated by all $u_\alpha -u_\beta$,
with $\alpha,\beta\in A$ such that $\alpha\sim\beta$. For every $u\in
U$ we denote by $\bar u$ its class in $U/W$.

Then $(\bar u_\alpha )_{\alpha\in B}$ is a basis for $U/W$.
\elm
\pf Let $U'\sbq U$ be the subspace generated by $u_\alpha$, with
$\alpha\in B$. Let $f:U\to U'$ be the linear function given by
$u_\alpha\mapsto u_\beta$, where $\beta$ is the unique element of $B$
such that $\alpha\sim\beta$. If $u_\alpha -u_\beta$, with
$\alpha\sim\beta$, is a generator of $W$ and $\gamma\in B$ such that
$\alpha\sim\gamma$ then we also have $\beta\sim\gamma$ and so
$f(u_\alpha )=f(u_\beta )=u_\gamma$. It follows that $u_\alpha
-u_\beta\in\ker f$ and so $W\sbq\ker f$. Therefore $f$ induces
a linear map $\bar f:U/W\to U'$, given by $\bar u_\alpha\mapsto
u_\beta$, where $\beta\in B$ such that $\alpha\sim\beta$.

We now define $g:U'\to U/W$, given by $u_\alpha\mapsto\bar u_\alpha$
$\forall\alpha\in B$. For every $\alpha\in B$ we have $\bar
fg(u_\alpha )=\bar f(\bar u_\alpha )=u_\alpha$. (We have $\alpha\in B$
and $\alpha\sim\alpha$.) Thus $\bar fg=1_{U'}$. If $\alpha\in A$ and
$\beta\in B$ such that $\alpha\sim\beta$ then $g\bar f(\bar u_\alpha
)=g(u_\beta )=\bar u_\beta =\bar u_\alpha$. (We have $\alpha\sim\beta$
so $u_\alpha -u_\beta\in W$ so $\bar u_\alpha =\bar u_\beta$.) Thus
$g\bar f=1_{U/W}$.

Thus $g:G'\to G/W$ is an isomorphism and $\bar f$ its
inverse. Since $(\bar u_\alpha )_{\alpha\in B}$ is the image with
respect to $g$ of the basis $(u_\alpha)_{\alpha\in B}$ of $G'$, it
will be a basis for $U/W$. \qed

\bpr Let $J=\{ (i_1,\ldots,i_n)\in I^n\mid i_1\leq\cdots\leq i_n\}$,
$J_1=\{ (i_1,\ldots,i_n)\in I^n\mid i_1<\cdots <i_n\}$ and $J_2=
J\setminus J_1$. Then
$$\{ v_{i_1}\odot\ldots\odot v_{i_n}\,\mid\, (i_1,\ldots,i_n)\in
J\}\cup \{ c(v_{i_1}\odot\ldots\odot v_{i_n})\,\mid\,
(i_1,\ldots,i_n)\in J_1\}$$
is a basis of $S'^n(V)$.
\epr
\pf We use Lemma 2.4 for $U=T^n(V)$, with the basis $(u_\alpha
)_{\alpha\in A}$, where $A=I^n$ and
$u_{i_1,\ldots,i_n}=v_{i_1}\otimes\cdots\otimes v_{i_n}$ $\forall
(i_1,\ldots,i_n)\in A$. The equivalence relation $\sim$ on $A$ is
given by $(i_1,\ldots,i_n)\sim (j_1,\ldots,j_n)$ if
$(j_1,\ldots,j_n)=(i_{\sigma (1)},\ldots,i_{\sigma (n)})$ for some
$\sigma\in A_n$ and $W\sbq U$ is generated by $u_\alpha -u_\beta$,
with $\alpha,\beta\in A$, $\alpha\sim\beta$.

We claim that $B=J\cup\{ (i_2,i_1,i_3\ldots,i_n)\mid
(i_1,\ldots,i_n)\in J_1\}$ is a set of representatives for
$A_{/\sim}$. First we show if $\alpha,\beta\in B$, $\alpha\neq\beta$,
then $\alpha\not\sim\beta$. We note that if $(i_1,\ldots,i_n)\sim
(j_1,\ldots,j_n)$ then then the sequences $i_1,\ldots,i_n$ and
$j_1,\ldots,j_n$ written in increasing order are the same. Therefore,
if $(i_1,\ldots,i_n),(j_1,\ldots,j_n)\in J$ with $(i_1,\ldots,i_n)\neq
(j_1,\ldots,j_n)$ then $(i_1,\ldots,i_n)$ or
$(i_2,i_1,i_3,\ldots,i_n)$ cannot be in the relation $\sim$ with
$(j_1,\ldots,j_n)$ or $(j_2,j_1,j_3,\ldots,j_n)$. This proves that if
$\alpha,\beta\in B$, with $\alpha\neq\beta$ then $\alpha\not\sim\beta$
unless $\alpha =(i_1,\ldots,i_n)$ and $\beta
=(i_2,i_1,i_3,\ldots,i_n)$ (or viceversa) for some
$(i_1,\ldots,i_n)\in J_1$, i.e. with $i_1<\cdots <i_n$. But in this
case the only $\sigma\in S_n$ such that
$(i_2,i_1,i_3,\ldots,i_n)=(i_{\sigma (1)},\ldots,i_{\sigma (n)})$ is
$\sigma =(1,2)$, which is odd. Hence, again, $\alpha\not\sim\beta$.

Next we prove that if $\alpha =(i_1,\ldots,i_n)\in A$ then there is
$\beta\in B$ with $\alpha\sim\beta$. We write the sequnce
$i_1,\ldots,i_n$ in increasing order as $j_1,\ldots,j_n$. Then
$(j_1,\ldots,j_n)\in J\sbq B$ and we have $(j_1,\ldots,j_n)=(i_{\sigma
(1)},\ldots,i_{\sigma (n)})$ for some $\sigma\in S_n$. If $\sigma\in
A_n$ then $(i_1,\ldots,i_n)\sim (j_1,\ldots,j_n)$ so we may take
$\beta =(j_1,\ldots,j_n)$. Suppose now that $\sigma\in S_n\setminus
A_n$. If $(j_1,\ldots,j_n)\in J_2$ then $j_k=j_{k+1}$ for some $1\leq
k\leq n-1$. Hence if $\tau =(k,k+1)\in S_n$ then
$(j_1,\ldots,j_n)=(j_{\tau (1)},\ldots,j_{\tau (n)})=(i_{\sigma\tau
(1)},\ldots,i_{\sigma\tau (n)})$. Since $\sigma,\tau\in S_n\setminus
A_n$ we have $\sigma\tau\in A_n$ so $(i_1,\ldots,i_n)\sim
(j_1,\ldots,j_n)$. So again we may take $\beta =(j_1,\ldots,j_n)$. If
$(j_1,\ldots,j_n)\in J_1$ the we also have $(j_2,j_1,j_3\ldots,j_n)\in
B$. If $\tau =(1,2)\in S_n$ then $(j_2,j_1,j_3,\ldots,j_n)=(j_{\tau
(1)},\ldots,j_{\tau (n)})=(j_{\sigma\tau (1)},\ldots,i_{\sigma\tau
(n)})$. Since $\sigma,\tau\in S_n\setminus A_n$ we have
$\sigma\tau\in A_n$ so $(i_1,\ldots,i_n)\sim
(j_2,j_1,j_3,\ldots,j_n)$. So this time we may take $\beta
=(j_2,j_1,j_3,\ldots,j_n)$.

The subspace $W_{S'}^n(V)$ of $T^n(V)$ is generated by $f_\sigma
(x_1,\ldots,x_n)$ with $x_1,\ldots,x_n$ and $\sigma\in A_n$, where
$f_\sigma :V^n\to T^n(V)$ is given by $(x_1,\ldots,x_n)\mapsto
x_1\otimes\cdots\otimes x_n-x_{\sigma (1)}\otimes\cdots\otimes
x_{\sigma(n)}$. But $f_\sigma$ is multilinear so we may restrict
ourselves to the case when $x_1,\ldots,x_n$ belong to the basis $v_i$,
with $i\in I$, of $V$. Then $W_{S'}^n(V)$ is generated by $f_\sigma
(v_{i_1},\ldots,v_{i_n})=v_{i_1}\otimes\cdots\otimes
v_{i_n}-v_{i_{\sigma (1)}}\otimes\cdots\otimes
v_{i_{\sigma(n)}}=u_{i_1,\ldots,i_n}-u_{i_{\sigma
(1)},\ldots,i_{\sigma (n)}}$, with $i_1,\ldots,i_n\in I$ and
$\sigma\in A_n$. Equivalently, $W_{S'}^n(V)$ is generated by $u_\alpha
-u_\beta$, with $\alpha,\beta\in A$, $\alpha\sim\beta$,
i.e. $W_{S'}^n(V)=W$. It follows that
$U/W=T^n(V)/W_{S'}^n(V)=S'^n(V)$. Also if $u=x_1\otimes\cdots\otimes
x_n\in U=T^n(V)$ then its class in $U/W=S'^n(V)$ is $\bar
u=x_1\odot\cdots\odot x_n$. In particular, $\bar
u_{i_1,\ldots,i_n}=v_{i_1}\odot\cdots\odot v_{i_n}$.

By Lemma 2.4, $\bar u_\alpha$ with $\alpha\in B$ are a basis of
$U/W=S'^n(V)$. If $\alpha =(i_1,\ldots,i_n)\in J$ then $\bar u_\alpha
=v_{i_1}\odot\cdots\odot v_{i_n}$. If $\alpha =(i_2,i_1,i_3,\ldots,i_n)$
for some $(i_1,\ldots,i_n)\in J_1$ then $\bar u_\alpha
=v_{i_2}\odot v_{i_1}\odot v_{i_3}\odot\cdots\odot
v_{i_n}=c(v_{i_1}\odot\cdots\odot v_{i_n})$. Hence the
conclusion. \qed

\btm For $n\geq 2$ we have the exact sequence
$$0\to\Lambda^n(V)\xrightarrow{\rho_{\Lambda^n,S'^n}}
S'^n(V)\xrightarrow{\rho_{S'^n,S^n}} S^n(V)\to 0,$$
where $\rho_{\Lambda^n,S'^n}$ is given by $x_1\wedge\cdots\wedge
x_n\mapsto x_1\odot\cdots\odot x_n- c(x_1\odot\cdots\odot x_n)$.
\etm
\pf We already know that $\rho_{S'^n,S^n}$ is surjective.

The map $(x_1,\ldots,x_n)\mapsto x_1\odot\cdots\odot
x_n-c(x_1\odot\cdots\odot x_n)$ is linear in each variable and
anti-symmetric. (If $x_i=x_j$ for some $i\neq j$ then, by Proposition
2.3(iv), we have $c(x_1\odot\cdots\odot x_n)=x_1\odot\cdots\odot
x_n$.) Hence the map $\rho_{\Lambda^n,S'^n}$, given by $x_1\wedge\cdots\wedge
x_n\mapsto x_1\odot\cdots\odot x_n-c(x_1\odot\cdots\odot x_n)$, is
well defined.

We prove the injectivity of $\rho_{\Lambda^n,S'^n}$. We use the
notations of Propostion 2.5. The set $\{ v_{i_1}\wedge\cdots\wedge
v_{i_n}\,\mid\, (i_1,\ldots,i_n)\in J_1\}$ is a basis of
$\Lambda^n(V)$. Let $\alpha\in\ker\rho_{\Lambda^n,S'^n}$. We write
$\alpha =\sum_{(i_1,\ldots,i_n)\in J_1}
a_{i_1,\ldots,i_n}v_{i_1}\wedge\cdots\wedge v_{i_n}$, with
$a_{i_1,\ldots,i_n}\in F$. Then
$$0=\rho_{\Lambda^n,S'^n}(\alpha )=\sum_{(i_1,\ldots,i_k)\in
J_1}a_{i_1,\ldots,i_k}(v_{i_1}\odot\cdots\odot
v_{i_n}-c(v_{i_1}\odot\cdots\odot v_{i_k})).$$
Since $v_{i_1}\odot\cdots\odot v_{i_n}$ and
$c(v_{i_1}\odot\cdots\odot v_{i_n})$, with $(i_1,\ldots,i_n)\in J_1$, are
part of the basis of $S'^n(V)$ form Proposition 2.5, this implies that
$a_{i_1,\ldots,i_n}=0$ $\forall (i_1,\ldots,i_n)\in J_1$ so
$\alpha =0$. Hence $\rho_{\Lambda^n,S'^n}$ is injective.

For the exactness in the second term we use the formulas
$\rho_{S'^n,S^n}(x_1\odot\cdots\odot
x_n)=\rho_{S'^n,S^n}c(x_1\odot\cdots\odot x_n)=x_1\cdots x_n$. (See
Propostion 2.3 (ii).)

The map $\rho_{S'^n,S^n}\rho_{\Lambda^n,S'^n}$ is given by
$x_1\wedge\cdots\wedge
x_k\mapsto\rho_{S'^n,S^n}(x_1\odot\cdots\odot
x_n-\\ c(x_1\odot\cdots\odot x_n))=x_1\cdots x_n-x_1\cdots x_n=0$. Hence
$\rho_{S'^n,S^n}\rho_{\Lambda^n,S'^n}=0$ so
$\ker\rho_{S'^n,S^n}\spq\Ima\rho_{\Lambda^n,S'^n}$.

For the reverse inclusion let $\alpha\in\ker\rho_{S'^n,S^n}$. By
Proposition 2.5, $\alpha$ writes as
$$\alpha =\sum_{(i_1,\ldots,i_n)\in
J}a_{i_1,\ldots,i_n}v_{i_1}\odot\cdots\odot
v_{i_n}+\sum_{(i_1,\ldots,i_n)\in
J_1}b_{i_1,\ldots,i_n}c(v_{i_1}\odot\cdots\odot v_{i_n}),$$
where $a_{i_1,\ldots,i_n},b_{i_1,\ldots,i_n}\in F$. Then we have
$$\begin{aligned}
0&=\rho_{S'^n,S^n}(\alpha )=\sum_{(i_1,\ldots,i_n)\in
J}a_{i_1,\ldots,i_n}v_{i_1}\cdots v_{i_n}+\sum_{(i_1,\ldots,i_n)\in
J_1}b_{i_1,\ldots,i_n}v_{i_1}\cdots v_{i_n}\\
&=\sum_{(i_1,\ldots,i_n)\in
J_1}(a_{i_1,\ldots,i_n}+b_{i_1,\ldots,i_n})v_{i_1}\cdots
v_{i_n}+\sum_{(i_1,\ldots,i_n)\in J_2}a_{i_1,\ldots,i_n}v_{i_1}\cdots
v_{i_n}.
\end{aligned}$$
Since $v_{i_1}\cdots v_{i_n}$ with $(i_1,\ldots,i_n)\in J=J_1\sqcup J_2$
are a basis of $S^n(V)$, we get $a_{i_1,\ldots,i_n}=0$
$\forall (i_1,\ldots,i_n)\in J_2$ and
$a_{i_1,\ldots,i_n}+b_{i_1,\ldots,i_n}=0$, so
$b_{i_1,\ldots,i_n}=-a_{i_1,\ldots,i_n}$, $\forall (i_1,\ldots,i_n)\in
J_2$. It follows that
$$\begin{aligned}
\alpha &=\sum_{(i_1,\ldots,i_n)\in
J_1}a_{i_1,\ldots,i_n}v_{i_1}\odot\cdots\odot
v_{i_n}+\sum_{(i_1,\ldots,i_n)\in 
J_1}-a_{i_1,\ldots,i_n}c(v_{i_1}\odot\cdots\odot
v_{i_n})\\
&=\sum_{(i_1,\ldots,i_n)\in 
J_1}a_{i_1,\ldots,i_n}(v_{i_1}\odot\cdots\odot v_{i_n}-c(v_{i_1}\odot\cdots\odot
v_{i_n}))=\rho_{\Lambda^n,S'^3}(\beta ),\end{aligned}$$
where $\beta=\sum_{(i_1,\ldots,i_n)\in
J_1}a_{i_1,\ldots,i_n}v_{i_1}\wedge\cdots\wedge v_{i_n}$. Thus
$\alpha\in\Ima\rho_{\Lambda^n,S'^3}$.\qed

Recall that if $n\leq 2$ then $S'^n(V)=T^n(V)$.

If $n=0,1$ then
$S^n(V)=S'^n(V)=T^n(V)$ and $\rho_{S'^n,S^n}$ is the identity map so
we have the short exact sequence $0\to 0\to
S'^n(V)\xrightarrow{\rho_{S'^n,S^n}}S^n(V)\to 0$. By putting together these
two sequences with those for $n\geq 2$ from Theorem 2.6, we get:

\bco We have a short exact sequence,
$$0\to\Lambda^{\geq 2}(V)\xrightarrow{\rho_{\Lambda^{\geq 2},S'}}S'(V)
\xrightarrow{\rho_{S',S}}S(V)\to 0.$$
\eco

{\bf Remark} The maps $\rho_{\Lambda^2,S'^2}$ and $\rho_{S'^2,S^2}$
are given by $x\wedge y\mapsto x\odot y-y\odot x$ and $x\odot y\mapsto
xy$, respectively. But when identify $S'^2(V)$ with $T^2(V)$ they
write as $x\wedge y\mapsto x\otimes y-y\otimes x=[x,y]$ and $x\otimes
y\mapsto xy$ so they coincide with $\rho_{\Lambda^2,T^2}$ and
$\rho_{T^2,S^2}$. Therefore the short exact sequences form Theorem
2.6, in the case $n=2$, and from Proposition 1.1 are the same.
\bigskip

{\bf References}
\medskip

\noindent [KT] Christian Kassel and Vladimir Turaev. ``Braid groups'', volume
247 of Graduate Texts in Mathematics. Springer, New York, 2008. With
the graphical assistance of Olivier Dodane.

\end{document}